\documentclass[11pt,a4paper]{amsart}
\usepackage{amscd,amssymb,graphicx}
\usepackage{amsfonts,amsmath,amssymb,epsf,epsfig}
\usepackage{xypic}
\xyoption{all}
\input{xypic}
\newdir{ >}{{}*!/-8pt/\dir{>}}
\def\N{\mathbb{N}}

\def\Z{\mathbb{Z}}

\def\fg{{\frak g}}

\def\fm{{\frak m}}

\def\pr{\noindent $\bf{Proof.}$\quad}     
\def\fin{\hfill$\square$\\}           
\newtheorem{theo}{Theorem}

\newtheorem{prop}{Proposition}
\newtheorem{cor}{Corollary}

\newtheorem{lem}{Lemma}
\author{Alice Fialowski}
\address{E\"otv\"os Lor\'and University\\ Budapest, Hungary}
\email{fialowsk@cs.elte.hu}

\author{Friedrich Wagemann}
\address{Universit\'e de Nantes\\ Nantes, France}
\email{wagemann@math.univ-nantes.fr}

\begin{document}
\title[Cohomology and deformations of \mbox{$\frak{m}_0$}] 
{Cohomology and deformations of the infinite dimensional filiform Lie 
algebra \mbox{$\frak{m}_0$}}

   
\keywords {Filiform Lie algebra, cohomology, deformation. Massey product }
\subjclass[2000]{ $17$B$65$, $17$B$56$, $58$H$15$}


\begin{abstract}
Denote $\fm_0$ the infinite dimensional $\N$-graded Lie algebra defined by 
basis $e_i$ $i\geq 1$ and relations $[e_1,e_i]=e_{i+1}$ for all $i\geq 2$.
We compute in this article the bracket structure on $H^1(\fm_0,\fm_0)$,
$H^2(\fm_0,\fm_0)$ and in relation to this, we establish that there are only
finitely many true deformations of $\fm_0$ in each nonpositive weight by constructing
them explicitely. It turns out that in weight 0 one gets exactly the other two 
filiform Lie algebras.
\end{abstract}
\maketitle

\section*{Introduction}

Recall the classification of infinite dimensional $\N$-graded Lie algebras 
$\fg=\bigoplus_{i=1}^{\infty}\fg_i$ with one-dimensional homogeneous components
$\fg_i$ and two generators over a field of characteristic zero. A. Fialowski 
showed in \cite{Fia1} that any Lie algebra of this type must be isomorphic to 
$\fm_0$, $\fm_2$ or $L_1$. We call these Lie algebras infinite dimensional 
filiform Lie algebras in analogy with the finite dimensional case where the 
name was coined by M. Vergne in \cite{Ver}. Here $\fm_0$ is given by generators
$e_i$, $i\geq 1$, and relations $[e_1,e_i]=e_{i+1}$ for all $i\geq 2$, $\fm_2$ 
with the same generators by relations $[e_1,e_i]=e_{i+1}$ for all $i\geq 2$, 
$[e_2,e_j]=e_{j+2}$ for all $j\geq 3$, and $L_1$ with the same generators is 
given by the relations $[e_i,e_j]=(j-i)e_{i+j}$ for all $i,j\geq 1$. $L_1$
appears as the positive part of the Witt algebra given by generators $e_i$ for
$i\in\Z$ with the same relations $[e_i,e_j]=(j-i)e_{i+j}$ for all $i,j\in\Z$.
The result was also obtained later by Shalev and Zelmanov in \cite{SZ}.

The cohomology with trivial coefficients of the Lie algebra $L_1$
was studied in \cite{Gont},
the adjoint cohomology in degrees $1$, $2$ and $3$ has been computed in 
\cite{Fia2} and also all of its non equivalent deformations were given. 
For the Lie algebra $\fm_0$, the cohomology with trivial coefficients 
has been studied in \cite{FialMill}, but neither the adjoint cohomology, nor 
related deformations have been computed so far. 
 The reason is probably that - as 
happens usually for solvable Lie algebras - the cohomology is huge and 
therefore meaningless. Our point of view is that there still remain interesting
features. We try to prove this in the present article by studying the adjoint
cohomology of $\fm_0$, while we reserve $\fm_2$ for a forthcoming paper.  

Indeed, it is true that the first and second adjoint cohomology of $\fm_0$ are 
infinite dimensional. The space $H^1(\fm_0,\fm_0)$ becomes already interesting 
when we split it up into homogeneous components $H^1_l(\fm_0,\fm_0)$ 
of weight $l\in\Z$, 
this latter space being finite dimensional for each $l\in\Z$. We compute the 
bracket structure on $H^1(\fm_0,\fm_0)$ in section $1$.

The space $H^2(\fm_0,\fm_0)$ is discussed in section $2$. This space is worse 
as it is infinite dimensional even in each weight 
separately. The interesting new feature here is that there are only finitely
many generators in each negative or zero weight which give rise to {\it true} 
deformations.
Given a generator of $H^2(\fm_0,\fm_0)$, i.e. an infinitesimal deformation,
corresponding to the linear term of a formal deformation, one can try to 
adjust higher order terms in order to have the Jacobi identity in the deformed
Lie algebra up to order $k$. If the Jacobi identity is satisfied to all orders,
we will call it a true (formal) deformation, see Fuchs' book \cite{Fuks} for
details on cohomology and [2] for deformations of Lie algebras.
  
In section 3 we discuss Massey products, in section 4 describe all true
deformations in negative weights. Section 5 deals with deformations in zero 
and positive weights.

As obstructions to infinitesimal deformations given by classes in 
$H^2(\fm_0,\fm_0)$ are expressed by Massey powers of these classes in 
$H^3(\fm_0,\fm_0)$, it is the vanishing of these Massey squares, cubes etc
which selects within the $H^2_l(\fm_0,\fm_0)$ of weight $l$ a finite number of
cohomology classes. The main result reads

\begin{theo}
The true deformations of $\fm_0$ are finitely generated in each weight
$l\leq 1$. More precisely, the space of unobstructed cohomology classes is 
in degree

\begin{itemize}
\item{$l\leq -3$} of dimension two,  
\item{$l= 0$} of dimension two,
\item{$l= -2$} of dimension three,
\end{itemize}

\noindent while there is no true deformation in weight $l=-1$. In weight
$l=0$, these are deformations to $\fm_2$ and $L_1$. In weight
$l=1$, there are exactly two true deformations, while in weight $l\geq 2$, 
there are at least two.
\end{theo}

We do not have more precise information about how many true deformations there 
are in positive weight, but there are always at least two. As a deformation in 
these weights is a true deformation
if and only if all of its Massey squares are zero (as cochains !), true 
deformations are determined by a countable infinite system of homogeneous 
quadratic equations in countably infinitely many variables. We didn't succeed 
in determining the space of solutions of this system.
 
We believe that the discussion of these examples of deformations are 
interesting as they go beyond the usual approach where the condition that    
$H^2(\fm_0,\fm_0)$ should be finite dimensional is the starting point for the
examination of deformations, namely the existence of a miniversal deformation
\cite{FiaFuc}.

Another attractive point of our study is the fact that in some cases the Massey
squares and cubes involved are not zero because of general reasons, but because
of the combinatorics of the relations. Thus the second adjoint cohomology of 
$\fm_0$ may serve as an example on which to study explicitely obstruction 
theory.

\noindent{\bf Acknowledgements:}\quad The work has been partially supported 
by the grants OTKA T034641 and T043034 and by the Erasmus program between
 E\"otv\"os Lor\'and University Budapest and Universit\'e Louis Pasteur Strasbourg.
Both authors are grateful to IHES where some of this work was accomplished,
to Yury Nikolayevsky for useful remarks and to Matthias Borer who helped us to
 get hold on the
weight $l>0$ case by MUPaD based computations.
   
\section{The space $H^1(\fm_0,\fm_0)$}

The Lie algebra $\fm_0$ is an $\N$-graded Lie algebra 
$\fm_0=\bigoplus_{i=1}^{\infty}(\fm_0)_i$ with $1$-di\-men\-sio\-nal graded 
components $(\fm_0)_i$ and generated in degree $1$ and $2$. Choosing a basis
$e_i$ of $(\fm_0)_i$, the only non-trivial brackets (up to skew-symmetry) read
$[e_1,e_i]=e_{i+1}$ for all $i$. 
We are computing in this section the first cohomology space $H^1(\fm_0,\fm_0)$
of $\fm_0$ with adjoint coefficients. As Lie algebra and module are graded, the
cohomology space splits up into homogeneous components, and we will always 
work with homogeneous cocycles $\omega(e_i)=a_ie_{i+l}$ for a scalar $a_i$ and
a given {\it weight} $l\in\Z$.

Concerning the cocycle identity $d\omega(e_j,e_i)=0$, let us first suppose 
that $j=1$ and $i>1$ (up to choosing the
symmetric case $j>1$ and $i=1$). In this case, it reads
$$\omega(e_{i+1})\,=\,[e_1,\omega(e_i)] - [e_i,\omega(e_1)]$$
or, putting in the expression of $\omega$, for all $l\geq 0$
$$a_{i+1}e_{i+l+1}\,=\,a_ie_{i+l+1} + \delta_{l,0}a_1e_{i+l+1}.$$
This means that for all $l\geq 0$, we must have
$$a_{i+1}\,=\,a_i + \delta_{l,0}\,\,\,a_1,$$
while for $l=-1$, we get the previous equation for $i\geq 3$ and $a_3=0$, 
for $l=-2$, we get the previous equation for $i\geq 4$ and $a_4=a_3=0$, 
and for $l\leq -3$, we get the previous equation for $i\geq -l+2$ and
$a_{-l+2}=a_{-l+1}=0$, while there is no equation for $i<-l$. 

The second situation where the cocycle identity has non-zero terms is
when $l\leq -1$, and $i$ and $j\geq 2$. In this case, there is 
only one non-zero term in the equation, and we get $a_i=0$ for $i+l=1$.

Now let us deduce the possible $1$-cocycles in different weights:

\noindent{\bf case 1}: $l\leq -1$

In case $l\leq -3$, the first identity means that all $a_i$ for 
$i\geq -l+2$ must be equal and $a_{-l+2}=a_{-l+1}=0$, therefore all
$a_i=0$ for $i\geq -l+1$, while there is no constraint on 
$a_1, a_2,\ldots, a_{-l}$. This is compatible with the second situation.

In case $l=-1$ and $l=-2$, the first constraint implies that all $a_i=0$ 
for $i\geq 3$, while there is no constraint on $a_1$ and $a_2$. 
The second identity is then already satisfied for $l\leq -2$, while for
$l=-1$, it implies $a_2=0$. 

But observe that the formula $\omega(e_i)=a_ie_{i+l}$ makes sense for 
$l\leq -1$ only if $i\geq -l+1$. Therefore all coefficients 
$a_1,\ldots,a_{-l}$ has to be set zero for $l\leq -1$. 

In conclusion, all cohomology is zero in weight $l\leq -1$. 

\noindent{\bf case 2}: $l\geq 1$

In this case, the cocycle identity means that all $a_i$ for $i\geq 2$ must be 
equal, while there is no constraint on $a_1$.

\noindent{\bf case 3}: $l=0$

In this case, the first identity means that all $a_i$ for $i\geq 3$ are 
determined by $a_1$ and $a_2$, while there is no constraint on $a_1$ and $a_2$.

Let us now examine the coboundaries: an element $x\in\fm_0$ determines a
$1$-coboundary by $\alpha_x(y):=[x,y]$ for all $y\in\fm_0$. In order to have
a homogeneous coboundary, we must take $x=e_i$ for some $i>0$; $\alpha_{e_i}$
is then homogeneous of weight $i$. Therefore we have:

\begin{itemize}
\item $dC^0_{l}(\fm_0,\fm_0)=0$ for $l\leq 0$
\item $dC^0_{l}(\fm_0,\fm_0)$ generated by $de_l=[e_l,-]$ for $l\geq 1$.
\end{itemize}

Observe that the coboundaries for $l\geq 2$ are non-zero only on $e_1$, thus
they can modify only the $a_1$-term of a cocycle. The coboundary for 
$l=1$ is zero on $e_1$ and non-zero and constant on all other $e_i$. It
thus kills the cocycle where all $a_i$ for $i\geq 2$ are equal.  

In conclusion, we have

\begin{theo}
$${\rm dim}\,\,H^1_l(\fm_0,\fm_0)\,=\,\left\{\begin{array}{ccc}
1 & {\rm for} & l\geq 1 \\
2 & {\rm for} & l=0 \\
0 & {\rm for} & l\leq -1 \end{array}\right.$$
\end{theo}

Let us now determine representatives of the non-zero cohomology classes:

In $H^1_0(\fm_0,\fm_0)$, we have the generators $\omega_1$ (corresponding to 
$a_1=1$ and $a_2=0$) and $\omega_2$ (corresponding to $a_1=0$ and $a_2=1$)
defined by:

$$\omega_1(e_k)\,=\,\left\{\begin{array}{ccc}
e_1 & {\rm for} & k=1 \\
0 & {\rm for} & k=2 \\
(k-2)e_k & {\rm for} & k\geq 3 \end{array}\right.$$

$$\omega_2(e_k)\,=\,\left\{\begin{array}{ccc}
0 & {\rm for} & k=1 \\
e_k & {\rm for} & k\geq 2 \end{array}\right.$$

In $H^1_l(\fm_0,\fm_0)$ for $l\geq 1$, we have two different kinds of 
cocycles: there is $\gamma$ for $l=1$, and $\alpha_l$ for $l\geq 2$:

$$\gamma(e_k)\,=\,\left\{\begin{array}{ccc}
ce_2 & {\rm for} & k=1 \\
0  & {\rm for} & k\geq 2 \end{array}\right.$$ 

$$\alpha_l(e_k)\,=\,\left\{\begin{array}{ccc}
0 & {\rm for} & k=1 \\
b_l e_{k+l} & {\rm for} & k\geq 2 \end{array}\right.$$ 

It is well known that $H^*(\fg,\fg)$ carries a graded Lie algebra structure for
any Lie algebra $\fg$, and that $H^1(\fg,\fg)$ forms a graded Lie subalgebra.
Let us compute this bracket structure on our generators:

Given $a\in C^p(\fg,\fg)$ and $b\in C^q(\fg,\fg)$, define
$$ab(x_1,\ldots,x_{p+q-1})\,=\,\sum_{\sigma\in{\rm Sh}_{p,q}}(-1)^{{\rm sgn}\,
\sigma}a(b(x_{i_1},\ldots,x_{i_q}),x_{j_1}\ldots,x_{j_{p-1}})$$
for $x_1,\ldots,x_{p+q-1}\in\fg$. The bracket is then defined by
$$[a,b]\,=\,ab-(-1)^{(p-1)(q-1)}ba.$$
It thus reads on $H^1(\fg,\fg)$ simply
$$[a,b](x)\,=\,a(b(x))-b(a(x)).$$

We compute

$$[\omega_1,\alpha_l](e_k)\,=\,\left\{\begin{array}{ccc}
0 & {\rm for} & k=1 \\
b_l(2+l-2)e_{2+l} & {\rm for} & k=2 \\
(k+l-2)b_le_{k+l}-(k-2)b_le_{k+l} & {\rm for} & k\geq 3 \end{array}\right.$$

Therefore $[\omega_1,\alpha_l]\,=\,l\alpha_l$.

$$[\omega_1,\gamma](e_k)\,=\,\omega_1(\delta_{k1}ce_2)-\gamma(\delta_{k1}e_1)$$

Therefore $[\omega_1,\gamma]\,=\,-\gamma$.

$$[\omega_2,\alpha_l](e_k)\,=\,\left\{\begin{array}{ccc}
0 & {\rm for} & k=1 \\
\omega_2(b_le_{k+l})-\alpha_l(e_k)=0 & {\rm for} & k\geq 2 \end{array}\right.$$

Therefore $[\omega_2,\alpha_l]\,=\,0$.

$$[\omega_2,\gamma](e_k)\,=\,\left\{\begin{array}{ccc}
\omega_2(ce_2)-0 & {\rm for} & k=1 \\
0-\gamma(e_k)=0 & {\rm for} & k\geq 2 \end{array}\right.$$ 

Therefore $[\omega_2,\gamma]\,=\,\gamma$.

$$[\alpha_l,\gamma](e_k)\,=\,\left\{\begin{array}{ccc}
\alpha_l(ce_2)-0 & {\rm for} & k=1 \\
0-\gamma(b_le_{k+l}))=0 & {\rm for} & k\geq 2 \end{array}\right.$$

This gives $[\alpha_l,\gamma]=\delta_{k1}cb_le_{2+l}$. This is a cocycle in 
weight $l+1$, $l\geq 2$, but by the list of coboundaries in weight $\geq 2$, 
we see that it is actually a coboundary. Therefore we have 
$[\alpha_l,\gamma]=0$ in cohomology.

$$[\omega_1,\omega_2](e_k)\,=\,\left\{\begin{array}{ccc}
0-0=0 & {\rm for} & k=1 \\
0-0=0 & {\rm for} & k=2 \\
(k-2)e_k-(k-2)e_k=0 & {\rm for} & k\geq 3 \end{array}\right.$$

Therefore $[\omega_1,\omega_2]\,=\,0$. It is also rather clear that 
$[\alpha_l,\alpha_m]=0$.

In summary:

\begin{theo}
The bracket structure on $H^1(\fm_0,\fm_0)$ is described as follows:
the commuting weight zero generators $\omega_1$ and $\omega_2$ act on the 
trivial Lie algebra generated by $\gamma$ in weight $1$ and the $\alpha_l$ 
for weight $l\geq 2$ as grading 
elements, $\gamma$ has degree $-1$ w.r.t. $\omega_1$, degree $1$ w.r.t. 
$\omega_2$, while $\alpha_l$ has degree $l$ w.r.t. $\omega_1$ and  degree $0$ 
w.r.t. $\omega_2$.
\end{theo}  

\section{The space $H^2(\fm_0,\fm_0)$}

Let us first compute $H^2(\fm_0,\fm_0)$. We work with homogeneous
cocycles $\omega(e_i,e_j)=a_{ij}e_{i+j+l}$ for a fixed weight $l\in\Z$, and for
$i,j\geq 1$, $i\not= j$.\\

\noindent{\bf 2.0} Observe that for weights $l\leq -3$, there are forbidden 
coefficients $a_{i,j}$, because they show up in front of $e_{i+j+l}$ with
$i+j+l\leq 0$. For example in $l=-3$, $a_{1,2}$ must be set to zero, in 
weight $l=-4$, $a_{1,2}$, $a_{1,3}$ must be set to zero, and so on.\\

\noindent{\bf 2.1} The cocycle identity reads 
\begin{eqnarray*}  
d\omega(e_i,e_j,e_k) &=& \omega([e_i,e_j],e_k) + \omega([e_j,e_k],e_i) + 
\omega([e_k,e_i],e_j) \\
&& - [e_i,\omega(e_j,e_k)] - [e_j,\omega(e_k,e_i)] -
[e_k,\omega(e_i,e_j)] =0.
\end{eqnarray*}
Let us first suppose that one index is equal to $1$. 
The identity reads then
\begin{equation}     \label{*}
(a_{i+1,j} + a_{i,j+1})e_{i+j+l+1} \,=\, (a_{i,j} 
- \delta_{j+l,0}\,\,a_{j,1} - \delta_{i+l,0}\,\,a_{1,i})e_{i+j+l+1},
\end{equation}
where $i,j\geq 2$, $i\not= j$. This identity makes only sense for 
$i+j+l\geq 2$, because in the above equation, the 
$a_{i,j}$ term shows up in front of the bracket of $e_1$ with $e_{i+j+l}$.
It is therefore not valid uniformly for all $i,j$ starting from $l\leq -4$. 

For $i+j+l<0$, there is no equation, while there is a special equation 
for $i+j+l=0,1$, namely 
$$a_{i+1,j} + a_{i,j+1} = 0.$$
Note that by {\bf 2.0}, coefficients $a_{i,j}$ with $i+j+l\leq 0$ are set 
to zero.\\
 
\noindent{\bf 2.2} The cocycle identity for $e_i,e_j,e_k$ for $i,j,k\geq 2$, 
$i\not=j$, 
$i\not=k$, and $j\not=k$, gives a non-zero factor only if $j+k+l=1$ or 
$i+j+l=1$ or $i+k+l=1$ (thus for $l\leq -4$). One can always arrange that
only one factor is possibly non-zero for given $i$, $j$ with $i+j+l=1$ 
(by choosing $k=\max(i+1,j+1)$, for example). Thus for weight $l\leq -4$, 
the coefficient $a_{i,j}$ with $i+j=-l+1$, $i\not=j$, $i,j\geq 2$, must be 
zero (which is compatible with the special equation !).\\

\noindent{\bf 2.3} Let us now consider coboundaries: expressing that 
$\omega\in Z^2_l(\fm_0,\fm_0)$
is a coboundary $\omega=d\alpha$ for some $1$-cochain 
$\alpha\in C^1_l(\fm_0,\fm_0)$, 
$\alpha(e_i)=\alpha_ie_{i+l}$ for all $i\geq 1$, gives by evaluation on $e_i$ 
and $e_j$:
$$a_{i,j}e_{i+j+l}\,=\,\alpha([e_i,e_j])-[e_i,\alpha(e_j)]+[e_j,\alpha(e_i)].$$
This equation makes sense only for $i+j+l\geq 1$ as all terms are multiples of 
$e_{i+j+l}$. Let us first take one index to be $1$, then we get
$$a_{1,i}\,=\,\alpha_{i+1} - \alpha_i + \alpha_1\,\,\delta_{l,0}$$
for all $i\geq \max(2,-l+2)$, because $\alpha_i$ appears in front of the 
bracket of $e_1$ with $e_{l+i}$. Thus all $a_{1,i}$, $i\geq \max(2,-l+2)$, 
can be taken to be zero by adding a coboundary.
For $i=-l,-l+1\geq 2$, we have the special equation 
$$a_{1,i}\,=\,\alpha_{i+1}.$$
It is now clear that, up to a coboundary, we may suppose 
for any $l$ that the last two terms in equation (\ref{*}) are zero.
Observe that the non-coboundary terms $a_{1,i}$
for $l\leq -3$ in a general cocycle, namely the terms with 
$i=2,\ldots,-l-1$, must be set to zero by {\bf 2.0}.\\

\noindent{\bf 2.4} For weight $l\leq -1$, we have additional coboundaries: 
indeed, there is 
a non-zero term in the coboundary equation for $e_i$, $e_j$, $i,j\geq 2$,
$i\not= j$ yielding
$$a_{-l+1,j}\,=\, -\alpha_{-l+1}$$
for all $j\geq 2$, $j\not=-l+1$. Be aware that the coefficient $\alpha_{-l+1}$
of the coboundary $d\alpha$ is linked to $a_{1,-l}$ by the equation 
$\alpha_{-l+1}=a_{1,-l}$ (cf {\bf 2.3}). Thus we cannot choose at the same
time to render $a_{1,-l}=0$ and $a_{-l+1,j}=0$ in weight $l\leq -1$ by
addition of a coboundary, we can impose only one of these conditions. 
This means for example that the cocycle given 
by coefficients $a_{i,j}$ with $a_{2,j}=1$ for all $j\geq 3$ and $a_{i,j}=0$ 
for all $i,j\not= 2$ (``the $2$-family'', cf {\bf 2.5}) is a coboundary in 
weight $l=-1$. Here $a_{1,1}=0$ and $\alpha_2$ are not linked.
More generally, the cochain given by coefficients $a_{i,j}$ 
with $a_{m+1,j}=1$ for all $j\geq m+2$, $a_{i,j}=0$ for all other $i,j>m+1$ 
(unless those which must be non-zero in order to respect antisymmetry) is 
cohomologuous to the cocycle consisting of the only non-trivial coefficient 
$a_{1,m}=1$ in weight $l=-m\leq -2$.\\  

Let us now reconsider the equations (\ref{*}) in the stable range, i.e. with
$i$ and $j$ such that $i,j\geq 2$, $i\not= j$, and $i+j\geq -l+2$:
\begin{equation}  \label{***}
a_{i+1,j} + a_{i,j+1} \,=\, a_{ij}
\end{equation}
We will adopt two different points of view on this system of equations:\\

\noindent{\bf 2.5 First point of view:}
Call the equations $a_{2,3}=a_{2,4}$, $a_{3,4}=a_{3,5}$, $a_{4,5}=a_{4,6}$,
 $\ldots$ , {\it diagonal equations}, and the terms involved 
{\it diagonal terms}. The prescription 
$$a_{i,i+1}=a_{i,i+2}=\left\{\begin{array}{ccc} 1 & {\rm for} & i=k \\
                                                0 & {\rm for} & i\not= k
\end{array}\right.$$
specifies uniquely (unicity is shown by induction) a solution to this system,
called the $k$th {\it family} or $k$-{\it series}. For the $k$th family, all
$a_{r,s}$ with $r>k$ are zero, $a_{r,s}=1$ for $r=k$, $a_{r,s}$ is linear
in $k$ for $r=k-1$ (and $s$ sufficiently big), $a_{r,s}$ is quadratic
in $k$ for $r=k-2$ (and $s$ sufficiently big), and so on.

Let us consider some examples, while we refer to section {\bf 4.5} and 
{\bf 5.1} for 
more informations; in the following expressions, all coefficients involving an 
index $1$ are set to zero, and the first non-zero column starting from the RHS
(i.e. the non-zero elements of the column $\{a_{m,k}\}_{k\geq m+1}$) 
is normalized to $1$.\\

{\bf The $2$ family:}\quad $a_{2,k}=1$, $a_{j,k}= 0$ 
for all $j,k\geq 3$.\\

{\bf The $3$ family:}\quad $a_{3,k}=1$, $a_{j,k}= 0$ 
for all $j,k\geq 4$, and $a_{2,3}=a_{2,4}=0$, $a_{2,k}= -(k-4)$ for 
all $k\geq 5$.\\

{\bf The $4$ family:}\quad $a_{4,k}=1$, $a_{j,k}= 0$ 
for all $j,k\geq 5$, $a_{3,4}=a_{3,5}=0$, $a_{3,k}= -(k-5)$ for 
all $k\geq 6$, and $a_{2,3}=a_{2,4}=a_{2,5}=a_{2,6}=0$, 
$a_{2,k}= \frac{(k-5)(k-6)}{2}$ for all $k\geq 7$ (even for all $k\geq 5$).\\

{\bf The $5$ family:}\quad  $a_{5,k}=1$, $a_{j,k}= 0$ 
for all $j,k\geq 6$, $a_{4,5}=a_{4,6}=0$, $a_{4,k}= -(k-6)$ for 
all $k\geq 7$, $a_{3,4}=a_{3,5}=a_{3,6}=a_{3,7}=0$, 
$a_{3,k}= \frac{(k-6)(k-7)}{2}$ for all $k\geq 8$ (even for all $k\geq 6$),
and $a_{2,3}=a_{2,4}=a_{2,5}=a_{2,6}=a_{2,7}=a_{2,8}=0$, 
$a_{2,k}= -\frac{(k-6)(k-7)(k-8)}{3!}$ for all $k\geq 9$ (even for all 
$k\geq 6$).\\

\noindent{\bf 2.6 Second point of view:} 
One can specify $a_{2,n}$ for all $n$, in such a way that the diagonal 
equations are satisfied. This implies that by choosing pairs 
$(a_{2,3},a_{2,4})$, $(a_{2,5},a_{2,6})$, $(a_{2,7},a_{2,8})$, and so on, the
first member is free, while the second member is determined by the 
corresponding diagonal equation.

Indeed, in $(a_{2,3},a_{2,4})$, $a_{2,4}$ is determined by $a_{2,3}=a_{2,4}$, 
in $(a_{2,5},a_{2,6})$, $a_{2,6}$ is determined by 
$a_{2,5}-a_{2,4}=a_{2,6}-a_{2,5}$ (which is just $a_{3,4}=a_{3,5}$),
in $(a_{2,7},a_{2,8})$, $a_{2,8}$ is determined by 
$((a_{2,7}-a_{2,6})-(a_{2,5}-a_{2,4}))=((a_{2,3}-a_{2,7})-(a_{2,7}-a_{2,6}))$
(which is just $a_{4,5}=a_{4,6}$).
All the other coefficients are then uniquely determined.

\noindent{\bf 2.7} In conclusion, it is clear that for each weight $l\in\Z$, 
there is a 
countably infinite number of independent $2$-cohomology classes. 
More precisely, in weight $l> -4$, the $k$ families with $k=2,3,\ldots$ 
represent
independent $2$-cohomology classes. In weight $l\leq -4$, {\bf 2.2}
shows that the $k$ family is contradictory for 
$$k < \left\{\begin{array}{ccc} 2 + \frac{-l-3}{2} & {\rm for} & 
l\,\,\,\,\,{\rm odd} \\ 2 + \frac{-l-2}{2} & {\rm for} & l\,\,\,\,\,
 {\rm even} \end{array}\right.$$
But there is still a countably infinite number of independent $2$-cohomology
classes in each weight.

\begin{theo}
$${\rm dim}\,\,H^2_l(\fm_0,\fm_0)\,=\,\infty$$
for each weight $l\in\Z$.
\end{theo}

\section{Massey products and deformations}

The $2$-cohomology is rather meaningless, as it is infinite dimensional even
in each weight separately. We ask now which of these homogeneous $2$-cocycles
gives rise to a deformation of $\fm_0$. A necessary condition is that the
class of the Massey square of the cocycle in question is zero. 
The first thing we will show is that for a large range of weights, 
even the condition that the Massey square is zero as a cochain 
is necessary and sufficient, and we will then determine all $2$-cocycles which 
have zero Massey square and are thus the infinitesimal part of a deformation 
of $\fm_0$ which is polynomial and of polynomial degree $1$. 
We will show in section $5$ that $\fm_0$ deforms (in a homogeneous way) 
to $\fm_2$
and to $L_1$, but to no other $\N$-graded  Lie algebra (non-isomorphic 
to $\fm_0$, $\fm_2$, and $L_1$). This is consistent with the classification
of $\N$-graded  Lie algebras with $1$-dimensional graded components, generated 
in degrees $1$ and $2$ \cite{Fia1}. 

Let $\omega$ be a $2$-cocycle, given as above by its coefficients $a_{i,j}$.
The {\it Massey square} of $\omega$ is by definition
\begin{equation}   \label{**}
M(a)_{ijk}\,=\,(a_{i,j}a_{i+j+l,k} + a_{j,k}a_{j+k+l,i} + a_{k,i}a_{k+i+l,j})
e_{i+j+k+2l}.
\end{equation} 
Observe that $M(a)_{ijk}=0$ if any two indices coincide.

\subsection{Massey squares and deformations}

\begin{prop}
Let $l\geq -1$.

If there exist $i,j,k$ with 
$M(a)_{ijk}\not= 0$, then it is a non-trivial $3$-cohomology class, and 
the $2$-cocycle $\omega$ is obstructed. Thus $\omega$ is obstructed if and
only if there exist $i,j,k$ with $M(a)_{ijk}\not= 0$.
\end{prop}

\pr Let $\alpha\in C^2(\fm_0,\fm_0)$ be a homogeneous $2$-cochain
with $\alpha(e_i,e_j)=b_{i,j}e_{i+j+m}$. Then 
$$d\alpha(e_i,e_j,e_k)\,=\,\alpha([e_i,e_j],e_k) - [e_i,\alpha(e_j,e_k)] + 
{\rm cycl.}.$$
Given a $2$-cocycle $\omega$ with non-zero Massey square $M(a)_{ijk}$, one
wants to find $\alpha$ which compensates $M(a)_{ijk}$, i.e. with 
$d\alpha(e_i,e_j,e_k)=M(a)_{ijk}$.
As $M(a)_{ijk}$ is of weight $2l$, one must have $m=2l$.

For $i,j,k\geq 2$, there is only one non-zero term in the coboundary equation
(cf {\bf 2.4}), and we have $d\alpha(e_i,e_j,e_k)=b_{i,j}e_{k+1}$ in case
$i+j+2l=1$ (the cases $j+k+2l=1$ or $i+k+2l=1$ are similar). Thus we can 
compensate all Massey squares $M(a)_{ijk}$ with $i+j+2l=1$, $j+k+2l=1$ or 
$i+k+2l=1$. As $i,j,k\geq 2$, the highest weight case appears for $l=-2$.

On the other hand, for $l\geq -1$ all $a_{1,i}$ can be taken to be zero 
by adding coboundaries (cf {\bf 2.3}). Thus $M(a)_{ijk}=0$ if one index is
equal to $1$, and the only squares to compensate are those with 
$i,j,k\geq 2$.\fin 

An interesting fact to note from the above proof is that the Massey squares 
that one can compensate by $3$-coboundaries are the $M(a)_{ijk}$ 
with $i+j+2l=1$, $j+k+2l=1$ or $i+k+2l=1$ for $i,j,k\geq 2$ in weight
$l\leq -2$. In the following, we will use the notation $M_{ijk}$ for the 
coefficient of $e_{i+j+k+2l}$ in the corresponding Massey square.

\section{Deformations in negative weights}

\subsection{True deformations in weight $-1$}

We now consider only {\it square zero cohomology} in weight $-1$, i.e. those 
classes in $H^2_{-1}(\fm_0,\fm_0)$ with Massey square equal to zero (not only
as a cohomology class, but as a cochain !). By the previous section, this 
determines all deformations of $\fm_0$ in weight $-1$.

First of all, the $2$-family (cf {\bf 2.5}) is a square zero $2$-cocycle, and
it is not contradictory in weight $-1$ (cf {\bf 2.7}). But the $2$-family is
actually a coboundary according to {\bf 2.4}.

Now assume that $a_{i,j}$ for all $i,j\geq 2$, $i\not= j$, defines a normalized
cocycle, i.e. $a_{1,s}=0$ for all $s$, which we may assume according to 
{\bf 2.3}. Observe that for $l=-1$, there is no special equation of type
$a_{1,-l}=\alpha_{-l+1}$. We assume further that all Massey squares 
$M_{ijk}$ are zero. Suppose that
$k$ is the first integer such that $a_{3,k}\not= 0$, $k\geq 4$. 

As all $3$ coefficients below $a_{3,k}$ are zero, equation (\ref{***}) shows
that the first non-zero $4$ coefficient is $a_{4,k-1}$, and that all $2$
coefficients are equal up to $a_{2,k}$, while $a_{2,k}\not=a_{2,k+1}$.

Denote $a_{2,3}=c$. We will establish a table for the coefficients in order
to examine the possible cases:

\begin{lem}
$a_{3,k}\not= 0$ implies $a_{3,k+1}\not= 0$.
\end{lem}

\pr If $a_{3,k+1}= 0$, then $M_{23k}=a_{2,3}a_{4,k}+a_{3,k}a_{k+2,2}$. But
by equation (\ref{***}), $a_{3,k}=a_{4,k}$ and 
$a_{2,k}-a_{3,k}=a_{2,k+1}=a_{2,k+2}$, and therefore
$M_{23k}=a_{2,3}a_{3,k}-a_{3,k}(a_{2,3}-a_{3,k})=a_{3,k}^2\not=0$.
This contradiction shows that $a_{3,k+1}\not= 0$.\fin

For $k=4$, $a_{3,k}=a_{3,k+1}$. Then $M_{234}=a_{3,4}(a_{2,4}-a_{2,6})$. If 
$a:=a_{3,4}\not=0$, then $a_{2,4}=a_{2,6}=:c$, and we have $a_{2,5}=c-a$ and 
thus $a-c=a+c$, implying $a=0$: contradiction. 

Let us now suppose $k>5$. 
Then $M_{34(k-1)}=a_{4,k-1}a_{k+2,3}=0$ implies $a_{k+2,3}=0$, because
$a_{4,k-1}=-a_{3,k}\not=0$ by equation (\ref{***}).
Consideration of $M_{23(k+1)}=0$ gives $a_{2,3}=a_{2,k+2}=a_{2,k+3}$. 
$M_{23(k+2)}=0$ gives $a_{4,k+2}=-a_{3,k+3}$. $M_{34k}=0$ implies that 
either $a_{3,k+3}=0$ or $a_{3,k+1}=2a_{3,k}$. But in this last case,
$a_{5,k-1}=-4a_{3,k}$, and $M_{35(k-1)}=0$ gives $a_{3,k+3}=0$ anyhow.
$M_{34(k+1)}=0$ gives $a_{3,k+4}=a_{4,k+3}$.

This gives the following table for the coefficients $a_{i,j}$ (where we 
used $i$ as the column index and $j$ as the row index, contrary to the usual 
convention) with $a_{3,k}=a$ and $a_{3,k+1}=b$: 

\vspace{1cm}
\hspace{2cm}
\begin{tabular}{|c||c|c|c|c|c|} \hline
    & 2   & 3  & 4   & 5       & 6        \\  \hline\hline
k-3 &  c  & 0  & 0   & 0       & -a       \\  \hline
k-2 &  c  & 0  & 0   & a       & 3a+b     \\  \hline
k-1 &  c  & 0  & -a  & -(2a+b) & -3a+b    \\  \hline
k   &  c  & a  & a-b & a-2b    & a-3b     \\  \hline
k+1 & c+b & b  & b   & b       & b        \\  \hline
k+2 &  c  & 0  & 0   & 0       & 0        \\  \hline
\end{tabular}
\vspace{1cm}

But the relation $a_{2,k}=c=a_{2,k+1}+a_{3,k}=c+b+a$ shows that $a=-b$.

Now, if $k$ is odd, we have $a_{3,k}=a$, $a_{5,k-2}=a$, $a_{7,k-4}=a$, and so 
on, until we reach the diagonal $a_{j,j}=a=0$. Thus in this case we have a 
contradiction.

But if $k$ is even, we will take the line of $a_{i,j}$ given by $i+j=4+k$, and
go to the diagonal: finally, we will also get $a=0$, i.e. a contradiction.

In conclusion, a non-zero $3$ coefficient for a square zero cocycle leads in 
weight $-1$ to a contradiction. But then all $2$ coefficients must be equal
by equations (\ref{***}), and this gives the $2$ family.

In conclusion, we have shown 

\begin{prop}
There is no non-trivial square-zero cohomology class in weight $l=-1$. 
In particular, there does not exist any non-trivial true cohomology 
class in weight $l=-1$.
\end{prop}  

\subsection{True deformations in weight $-2$}

We saw in the last section how to determine all cocycles, here in weight $-2$,
which lead to true deformations. But as proposition $1$ in section 
$3.1$ is not valid in weight $-2$, we cannot use the vanishing of all 
Massey squares to get restrictions on our cocycle, instead, we have to leave
out those which are coboundaries and could thus be compensated by higher
Massey products. 

Let $\omega$ be a cocycle given by its coefficients $a_{i,j}$. Note that we 
can still suppose $a_{1,j}=0$ for all $j\geq -l+1$ and all $j\leq -l-1$
(cf {\bf 2.0}, {\bf 2.3}, {\bf 2.4}). 
On the other hand,
we may suppose that we are in the stable range and by {\bf 2.4} that the 
first terms in the $3$ column (i.e. at least $a_{3,4}$, $a_{3,5}$) are zero
up to a coboundary. As we cannot assume simultaneously that $a_{3,4}$, 
$a_{3,5}$ and $a_{1,-l}=a_{1,2}$ are zero (cf {\bf 2.4}), we choose to allow
$a_{1,2}$ non-zero. 

When writing a Massey square $M_{ijk}$, we will now always suppose that the 
indices are ordered $i<j<k$. The Massey squares which may be compensated are 
those $M_{ijk}$ with $i+j+2l=1$, according to section {\bf 3.1}. This means in 
weight $l=-2$ that all $M_{23k}$ can be compensated, and that these are the 
only (ordered) ones. The other $M_{ijk}$ must be zero.    

We start now a case study in order to determine which possibilities there are 
for $\omega$, imposing that all (ordered) Massey squares $M_{ijk}$ with 
$i+j\not=5$ are zero.\\

\noindent\quad{\bf 1st case:}\quad Suppose $a_{2,3}=a_{2,4}=0$. Then 
$a_{2,5}=-a_{3,4}=-a_{3,5}$ and $a_{2,5}=a_{3,5}+a_{2,6}$ implying that
$a_{2,6}=-2a_{3,5}$. $M_{245}=a_{4,5}(a_{2,5}-a_{2,7})=0$.

{\bf case 1a:}\quad $a_{2,5}=0$, and thus $a_{3,4}=a_{3,5}=0$, $a_{2,6}=0$. 
Then either $a_{4,5}=0$ ($\Rightarrow$ $a_{3,6}=0$, $a_{2,7}=0$), or 
$a_{2,7}=0$ ($\Rightarrow$ $a_{3,6}=0$, $a_{4,5}=0$). In any case $a_{2,7}=0$, 
$a_{3,6}=0$ and $a_{4,5}=0$. 

Suppose now given $r$ ($r\geq 10$) such that $a_{i,j}=0$ for all $i+j\leq r$.
Then $M_{2ij}=a_{i,j}(a_{2,i}-a_{2,i+j-2}+a_{2,j})$ must be zero for $i\geq 4$.
Let us suppose $i<j$ (indices ordered !) and $i+j=r+1$. Then by hypothesis
$a_{2,i}=a_{2,j}=0$ and $M_{2jk}=-a_{i,j}a_{2,r-1}$. Thus either $a_{i,j}=0$,
or $a_{2,r-1}=0$. But these two elements are on a new diagonal (in the 
matrix of coefficients $a_{i,j}$), and all elements with lower indices are 
zero. By equation (\ref{***}) this implies that two (because approaching the 
diagonal, one jumps to the next diagonal by the diagonal equations (cf 
{\bf 2.5}) $a_{s,s+1}=a_{s,s+2}$) new diagonals are zero, and by induction,
all coefficients are zero in this case.

{\bf case 1b:}\quad $a_{2,5}=1$, and thus $a_{2,6}=2$, $a_{3,4}=a_{3,5}=-1$.
Now $M_{245}=0$ implies $a_{4,5}=0$ or $a_{2,7}=1$. But for $a_{2,7}=1$, we get
by repeated use of equation (\ref{***}) $a_{3,6}=1$, $a_{4,5}=-2$, 
$a_{4,6}=-2$, $a_{3,7}=3$ and finally $a_{2,8}=-2$. Relate this then to
$M_{246}=a_{4,6}(-a_{2,8}+a_{2,6})\not=0$, to conclude that $a_{2,7}\not=1$,
and therefore $a_{4,5}=0$.

This means that all $a_{i,j}$ with $i+j\leq 10$ are the same as for the 
$3$-family (cf {\bf 2.5}). Let us show the induction step in order to conclude 
that two more diagonals are like in the $3$-family. Indeed, suppose now
$a_{i,j}$ with $i+j\leq r$ like in the $3$-family, and take $j=r-3$. We have 
$M_{24j}=a_{4,j}(a_{2,j}-a_{2,j+2})=a_{4,r-3}(a_{2,r-3}-a_{2,r-1})=0$. 
The coefficient $a_{2,r-3}$ must be as in the $3$-family by hypothesis. We 
want to conclude that $a_{4,r-3}=0$ (as in the $3$-family), opening up two 
more diagonals. Therefore we show that $a_{2,r-3}\not=a_{2,r-1}$. 

Let us denote $a_{2,r-3}=t$. We have by hypothesis $a_{2,r-2}=t+1$, $a_{3,r-4}
=a_{3,r-3}=-1$, and $a_{k,s}=0$ for $k<s$, $k\geq 4$. Suppose 
$a_{2,r-3}=a_{2,r-1}$, and this will lead to $a_{4,r-2}\not=0$ while 
$a_{2,r-2}=t+1$ and $a_{2,r}\not=t+1$. More precisely, in the new diagonal
starting from $a_{2,r-1}=t$, we get $a_{3,r-2}=1$, $a_{4,r-3}=-2$, 
$a_{5,r-4}=2$, and then we always get $\pm 2$, because there are only zeroes 
one diagonal higher. By construction $r-1$ is odd, say $r-1=2k+1$. Doing
in this sense $k-1$ steps on the diagonal towards the diagonal transforms
$a_{2,r-1}=a_{2,2k+1}$ into $a_{k-1+2,k+2}=a_{k+1,k+2}$. But by the diagonal 
equation (cf {\bf 2.5}), $a_{k+1,k+2}=a_{k+1,k+3}$, and then we work back 
$k-3$ steps to get $\pm 2(k-3)$, and finally $a_{2,r}=-(\mp 2(k-3)+1)+t$. 
This is equal to $t+1$ only if $k=2$ or $k=4$. $k=2$ is already treated, and 
for $k=4$, one can check directly that $a_{2,8}\not=a_{2,10}$:

\vspace{1cm}
\hspace{3cm}
\begin{tabular}{|c|c|c|c|c} \hline
2   &  3   &  4   &  5   \\  \hline\hline
0   &      &      &     \\  \hline
0   &  -1  &      &     \\  \hline
1   &  -1  &  0   &    \\  \hline
2   &  -1  &  0   &  2  \\  \hline
3   &  -1  &  -2  &  2  \\  \hline
4   &   1  &  -4  &    \\  \hline
3   &   5  &      &    \\  \hline
-8  &      &      &    \\  \hline
\end{tabular}
\vspace{1cm}
  
In conclusion, the only non-zero cocycle (making zero the non-compensable 
Massey squares) compatible with case $1$ is the $3$-family.\\ 

\noindent\quad{\bf 2nd case:}\quad Here we can take $a_{2,3}=a_{2,4}=1$.
Recall that we choose to take the first terms in the $3$ column (i.e. at least 
$a_{3,4}$, $a_{3,5}$) to be zero (possibly by adding a coboundary).

{\bf case 2a:}\quad Suppose as a first subcase $a_{4,5}=a_{4,6}=0$. Then we 
have up to $a_{i,j}$ with $i+j=10$ the $2$-family. Set $a_{5,6}=a$. We get then
 $a_{2,9}=1-a$ and $a_{2,10}=1-4a$ by repeated use of equation (\ref{***}).
But $M_{247}=a_{4,7}(a_{2,4}-a_{2,9}+a_{2,7})=0$ implies $a=0$ or $a=-1$, 
while $M_{248}=a_{4,8}(a_{2,4}-a_{2,10}+a_{2,8})=0$ implies $a=0$ or $a=-
\frac{1}{4}$. In conclusion, $a=0$ and the $2$-family is reproduced one 
diagonal higher. Using $M_{24j}$ and $M_{24(j+1)}$, one can show in a 
similar way that the only solution here is the $2$-family. 

{\bf case 2b:}\quad Here $a_{2,3}=a_{2,4}=1$, $a_{3,4}=a_{3,5}=0$, but 
$a:=a_{4,5}=a_{4,6}\not=0$. By $M_{245}=a_{4,5}(a_{2,4}-a_{2,7}+a_{2,5})$, 
$M_{246}=a_{4,6}(a_{2,4}-a_{2,8}+a_{2,6})$ and
$M_{345}=a_{4,5}(a_{3,5}-a_{3,7})$, we get thus $a_{3,5}=a_{3,6}=0$, 
$a_{2,4}+a_{2,5}=a_{2,7}$ and $a_{2,4}+a_{2,6}=a_{2,8}$. But then
on the one hand $a_{2,3}=a_{2,4}=a_{2,5}=a_{2,6}=a_{2,7}$ by equation 
(\ref{***}), but on the other hand $a_{2,4}+a_{2,5}=a_{2,7}=2$, which is a 
contradiction.

As a conclusion of the case study, the only cocycles which can possibly give
true deformations are the $2$-and the $3$-family, but possibly with a 
non-zero $a_{1,2}$ coefficient.\\

The $2$-family (cf {\bf 2.5}) is a square zero $2$-cocycle in weight $l=-2$, 
it is not contradictory in weight $-2$ (cf {\bf 2.7}), and is thus one solution
here. 

The $3$-family $\omega$ has a non-zero Massey square, namely 
$M_{23j}=a_{2,j}-a_{2,j+1}=1$ for all $j\geq 4$. Let us show that the 
corresponding Massey cube is then zero, and thus that the $3$-family gives 
indeed rise to an true deformation in weight $-2$: 

We must write $M_{23k}$ as a coboundary. The cochain $\alpha(e_i,e_j)=b_{i,j}
e_{i+j-4}$ with $b_{i,j}=M_{23k}$ for $i=2$ and $j=3$, $b_{i,j}=M_{32k}$ for 
$i=3$ and $j=2$, and $b_{i,j}=0$ otherwise satisfies
$$d\alpha(e_i,e_j,e_k)\,=\,M_{23k}e_{k+1}.$$
We must then compute the Massey cube
$$N_{ijk}:=\alpha(\omega(e_i,e_j),e_k)+\omega(\alpha(e_i,e_j),e_k) + 
{\rm cycl.} \,=\,a_{i,j}b_{i+j-2,k} + b_{i,j}a_{i+j-2,k}\,\,+\,\,{\rm cycl.}.$$
But if $b_{i+j-2,k}\not=0$, then $i+j-2,k\in\{2,3\}$. The only possibly non 
zero term is thus $N_{23k}=M_{23k}a_{1,k}=0$ ($k\geq 4$ here).

Finally, let us show that we cannot get any information about $a_{1,2}$, 
neither by the cocycle equations, nor by the vanishing of the Massey squares.
This is clear for the cocycle equations. Let us show that we cannot deduce
$a_{1,2}=0$ from Massey squares which have to vanish. Indeed, when writing down
the Massey squares which involve $a_{1,2}$, the only possibly non-zero Massey 
squares $M_{ijk}$ involving $a_{1,2}$ (with ordered indices) have $i=1$. But 
then we have to have $j=2$ in order to involve $a_{1,2}$. One easily checks 
that $M_{12k}=0$.    

To summarize, we have the following

\begin{prop}
The $3$-family in weight $-2$ has a non-zero Massey square, but
its Massey cube is zero, and we get consequently a true deformation.
In weight $-2$, the $2$- and $3$-family with possibly a non-zero term $a_{1,2}$
define the only cohomology classes leading to true deformations. 
\end{prop}
 
\subsection{True deformations in weight $-3$}

We will determine all cocycles leading to true deformations in weight 
$-3$ once again by imposing on a general cocycle $\omega$ given by its 
coefficients $a_{i,j}$ for all $i,j\geq 2$, $i\not= j$, that all Massey 
squares which cannot possibly be compensated (cf section {\bf 3.1}) are zero. 
The squares which cannot serve to give conditions on the $a_{i,j}$ are 
those $M_{ijk}$ (with ordered indices $i<j<k$) with $i+j=7$, $i+k=7$ or 
$j+k=7$.

All $1$-coefficients other than $a_{1,2}$ and $a_{1,3}$ may be supposed to be 
zero by {\bf 2.3} (cf {\bf 2.4}). $a_{1,2}=0$ by {\bf 2.0}. We choose once 
again that the first $4$-coefficients (i.e. at least $a_{4,5}=a_{4,6}$) are 
zero, up to a coboundary, according to {\bf 2.5}, while not imposing anything
on $a_{1,3}$ (cf {\bf 2.4}).

Let us draw the table for the coefficients of $\omega$:

\vspace{1cm}
\hspace{3cm}
\begin{tabular}{|c|c|c|c|c} \hline
2    &  3   &  4   &  5   \\  \hline\hline
a    &      &      &     \\  \hline
a    &  b   &      &     \\  \hline
a-b  &  b   &  0   &     \\  \hline
a-2b &  b   &  0   &     \\  \hline
a-3b &  b   &      &     \\  \hline
a-4b &      &      &     \\  \hline
\end{tabular}
\vspace{1cm}

Now we write down the Massey squares that we may use: 
$M_{236}=a_{2,6}(a_{2,3}-a_{3,6}+a_{3,5})=a(a-2b)$, 
$M_{237}=a_{2,7}(a_{2,3}-a_{3,7}+a_{3,6})=a(a-3b)$,
$M_{246}=a_{2,4}a_{3,6}+a_{4,6}a_{7,2}+a_{6,2}a_{5,4}=ab$.

In conclusion, $a=0$. But then up to $i+j=10$, the $3$-family has built up.
Let us show by induction that the $3$-family is the only possible solution:

Suppose the $a_{i,j}$ up to $i+j=r$ for $r\geq 10$ are like in the $3$-family
(cf {\bf 2.5}). Consider the Massey square
$$M_{23k}\,=\,a_{2,3}a_{2,k}+a_{3,k}a_{k,2}+a_{k,2}a_{k-1,3}\,=\,a_{2,k}
(a_{2,3}-a_{3,k}+a_{3,k-1}).$$
We may use its vanishing to deduce restrictions on the $a_{i,j}$ as soon as
$k\geq 6$. For $k\geq 8$ and with $r=k+2$, $M_{23k}=0$ implies under the
induction hypothesis that $a_{3,k}=a_{3,k-1}$, and we have therefore 
transmitted the $3$-family to two more diagonals, showing the induction 
step. \\

In order to conclude, let us show that the $3$-family is of Massey square 
zero:\\

Recall that the $3$-family is defined by $a_{3,k}=a\not= 0$, $a_{j,k}= 0$ 
for all $j,k\geq 4$, and $a_{2,3}=a_{2,4}=0$, $a_{2,k}= -(k-4)a$ for 
all $k\geq 5$.

It is clear that $M(a)_{ijk}=0$ for all $i,j,k\geq 4$, by definition of the 
$3$-family. Suppose $i=3$ ($j=3$ or $k=3$ would be a symmetric case):
$$M(a)_{3jk}=a_{3,j}a_{j,k}+a_{j,k}a_{j+k-3,3}+a_{k,3}a_{k,j}
=a_{j,k}(a_{3,j}+a_{j+k-3,3}+a_{3,k}).$$
This last expression is zero if both $j$ and $k$ are greater or equal to 
$4$ (as then $a_{j,k}=0$), and also if one of them is equal to $2$, because 
in this case the term in parenthesis is zero. Suppose now that $i=2$.
$$M(a)_{2jk}=a_{2,j}a_{j-1,k}+a_{j,k}a_{j+k-3,2}+a_{k,2}a_{k-1,j}.$$
In case $j,k\geq 4$, this expression reduces to 
$a_{2,j}a_{j-1,k}+a_{k,2}a_{k-1,j}$ which is evidently zero if both $j$ and 
$k$ are greater or equal to $5$, and in case $j=4$, $a_{k-1,j}$ and $a_{2,j}$ 
are zero. It remains the case where $j=3$, but then we get
$a_{3,k}a_{k,2}+a_{k,2}a_{k-1,3}=0$.

Finally, let us show that the possibly non-zero coefficient $a_{1,3}$ cannot 
be shown to be zero using the vanishing of Massey squares. The only $M_{ijk}$
(with ordered indices) involving $a_{1,3}$ have $i=1$. 
$$M_{1jk}\,=\,a_{1,j}a_{j+1+l,k}+a_{j,k}a_{j+k+l,1}+a_{k,1}a_{k+1+l,j}.$$
Then for $j=3$, we get $M_{13k}=0$, and for $j+k=6$, we get also $M_{1jk}=0$,
and these are the only combinations (up to reordering) involving $a_{1,3}$. 

To summarize, we get the following

\begin{prop}
In weight $-3$, the $3$-family, with a possibly non-zero $a_{1,3}$ coefficient,
defines the only cohomology class leading to an true deformation.
\end{prop}

\subsection{True deformations in weight $-4$}

We will determine all cocycles leading to true deformations in weight 
$-4$ once again by imposing on a general cocycle $\omega$ given by its 
coefficients $a_{i,j}$ for all $i,j\geq 2$, $i\not= j$ that all Massey squares 
which cannot possibly be compensated (cf section {\bf 3.1}) are zero. 
The squares which cannot serve to give conditions on the $a_{i,j}$ are 
those $M_{ijk}$ (with ordered indices $i<j<k$) with $i+j=9$, $i+k=9$ or 
$j+k=9$. In weight $l=-4$, we have to be more careful with the conditions 
as we are not always in the stable range (cf {\bf 2.5}). 
For example, {\bf 2.1} implies 
here that $a_{2,4}=0$ (and we can not deduce here $a_{2,3}=a_{2,4}$).
But for $j>i\geq 3$, and for $i=2$ and $j\geq 4$ we still have
$$a_{i+1,j}+a_{i,j+1}\,=\,a_{i,j}.$$ 

But then {\bf 2.2} implies that $a_{2,3}=0$. All $1$-coefficients other than
$a_{1,2}$, $a_{1,3}$ and $a_{1,4}$ may be supposed to be zero by {\bf 2.3}, 
while $a_{1,2}=a_{1,3}=0$ follows from {\bf 2.0}. We choose once again 
according to {\bf 2.4} that the first $5$-coefficients (i.e.
at least $a_{5,6}=a_{5,7}$) are zero, up to a coboundary, while we do not 
impose anything on $a_{1,4}$, cf {\bf 2.4}.

Let us draw the table for the coefficients of $\omega$:

\vspace{1cm}
\hspace{3cm}
\begin{tabular}{|c|c|c|c|c} \hline
2    &  3   &  4   &  5   \\  \hline\hline
0    &      &      &     \\  \hline
0    &  b   &      &     \\  \hline
-b   &  b   &  c   &     \\  \hline
-2b  & b-c  &  c   &  0  \\  \hline
c-3b & b-2c &  c   &  0  \\  \hline
3c-4b& b-3c &  c   &     \\  \hline
6c-5b& b-4c &      &     \\  \hline
10c-6b&      &      &     \\  \hline
\end{tabular}
\vspace{1cm}

Now we write down the Massey squares that we may use: 
$M_{246}=a_{4,6}a_{6,2}=2cb$, 
$M_{238}=a_{3,8}a_{7,2}+a_{8,2}a_{6,3}=(b-3c)(3b-c)+(3c-4b)(b-c)=-b(b+3c)$,
$M_{256}=a_{2,5}a_{3,6}+a_{5,6}a_{7,2}+a_{6,2}a_{4,5}=b(-b+3c)$.

Now start a case study: either $b=0$, and in this case we want to show that 
the $4$-family is built up by induction. Indeed, we have 
$M_{24k}=a_{2,k}(a_{4,k-2}-a_{4,k})$ which must vanish as soon as $k\geq 8$. 
In this way we transmit the built up $4$-family to another two diagonals.
Or $c=0$ and in this case all coefficients are zero. The zero family is also 
easily shown to be built up from this initial stage. 

Let us show that the nullity of $a_{1,4}$ cannot be derived from the nullity 
of Massey squares. The Massey squares (with ordered indices) where $a_{1,4}$ 
shows up, have either $i=1$ or they are $M_{234}$. The latter is zero anyhow,
and the former are shown to be zero as for $l=-1,-2$ and $-3$. 

To summarize, we have the following

\begin{prop}
The only cohomology class leading to an true deformation in weight 
$l=-4$ is represented by the $4$-family, with a possibly non-zero coefficient
$a_{1,4}$.
\end{prop} 

The fact that the $4$-family has zero Massey square in weight $-4$
follows from Proposition $7$ in section {\bf 4.5}.   

\subsection{True deformations in weight $l$, $l\leq -5$}

We will show that in degree $l$, the $-l$ family is of Massey square zero, 
and that this is the only family for which all Massey squares which cannot be
compensated, are zero. Therefore, we will show that in weight $l$, $l\leq -5$,
the $-l=:m$ family, with a possibly non-zero coefficient $a_{1,m}$, is the 
only cocycle which leads to true deformations.
  
For this, we need the explicit expression of the non-zero
low degree coefficients of the $-l$ family. It is obvious from {\bf 2.5} how 
to deduce the expressions of the coefficients 
of the general $m:=-l$ family from those for the low degree families:\\

{\bf The $m$ family:}\quad  $a_{m,k}=1$, $a_{j,k}= 0$ 
for all $j,k\geq m+1$, $a_{m-1,m}=a_{m-1,m+1}=0$, $a_{m-1,k}= -(k-(m+1))$ for 
all $k\geq m+2$, $a_{m-2,m-1}=a_{m-2,m}=a_{m-2,m+1}=a_{m-2,m+2}=0$, 
$a_{m-2,k}= \frac{(k-(m+1))(k-(m+2))}{2}$ for all $k\geq m+3$,
and $a_{m-3,m-2}=a_{m-3,m-1}=a_{m-3,m}=a_{m-3,m+1}=a_{m-3,m+2}=a_{m-3,m+3}=0$, 
$a_{m-3,k}= -\frac{(k-(m+1))(k-(m+2))(k-(m+3))}{3!}$ for all $k\geq m+4$,
and so on. Observe that the general expression of a coefficient $a_{i,j}$
in this family

\begin{prop}
The $m$-family defines a $2$-cocycle in any weight.
\end{prop}

\pr We have to show that the $m$ family  
satisfies the requirements of sections {\bf 2.1} (i.e. equation (\ref{***});
observe that with the non-zero coefficients of the $m$-family, we are always 
in the stable range) and {\bf 2.2}. It is clear that the requirement of 
{\bf 2.2} is met.

For the equation (\ref{***}), take the general expression 
of the above coefficients 
$$a_{m-r,k}\,=\,\pm\frac{(k-(m+1))!}{r!(k-(m+r+1))!}$$
for all $k\geq m+r+1$, and all $r\leq m-2$; $\pm$ denotes an alternating sign 
with respect to the parity of $r$. Now
$$a_{m-r,k+1}+a_{m-(r-1),k}=\frac{(k+1-(m+1)-r)(k-(m+1))!}{r!(k-(m+r))!}=
a_{m-r,k}.$$
\fin

Let us show now that the Massey square of the $m$-family is zero (i.e. not only
the non-compensable Massey squares, but all). 

\begin{prop}
All Massey squares of the $m$-family are zero in weight $l=-m$.
\end{prop}

\pr Indeed, we have
$$M_{ijk}\,=\,a_{i,j}a_{i+j+l,k}+a_{j,k}a_{j+k+l,i}+a_{k,i}a_{k+i+l,j}.$$ 
We will always consider ordered Massey squares, i.e. $M_{ijk}$ with $i<j<k$,
and it will be enough to show that these are zero. 

\noindent{\bf First case:}\quad $i=m-r$, $j=m-p$, and $k\geq m+1$ with 
$p,r\geq 0$. These conditions imply $a_{i,j}=0$, and we get
\begin{eqnarray*}
M_{ijk}&=& a_{j,k}a_{j+k+l,i}+a_{k,i}a_{k+i+l,j} \\
&=& (-1)^{p+1}\frac{(k-(m+1))!}{p!(k-(m+p+1))!}
(-1)^{r}\frac{(k-p-(m+1))!}{r!(k-p-(m+p+1))!}  \\
&+&(-1)^{r+1}\frac{(k-(m+1))!}{r!(k-(m+r+1))!}
(-1)^{p+1}\frac{(k-r-(m+1))!}{p!(k-r-(m+p+1))!}
\end{eqnarray*}
Suppose now first that $k+i+l=k-r>j$. In this case we get by taking out 
common factors
$$\textstyle{
M_{ijk}=(-1)^{p+r+1}\frac{(k-(m+1))!}{r!p!(k-p-(m+r+1))!}
\left(\frac{(k-p-(m+1))!}{(k-(m+p+1))!}-\frac{(k-r-(m+1))!}{(k-(m+r+1))!}
\right)}$$
which is zero. On the other hand, in case $k+i+l=k-r=j$, we get $k=m+q$,
$i=m-r$, $j=m-p$. Then the only possibly non-zero term is 
$M_{ijk}= a_{j,k}a_{j+k+l,i}$, because $a_{k+i+l,j}=0$. But 
$a_{j,k}a_{j+k+l,i}=a_{m-p,m+q}a_{m+(q-p),m-(p+q)}$ and 
$a_{m+(q-p),m-(p+q)}=0$, because $a_{m-s,m+s}=0$, $a_{m-s,m+s+1}\not=0$ 
marks the last zero term in the $m$ family (when fixing $m-s$ and counting up 
the second index), but here $q-p\leq p+q$. It remains the third subcase where
$k-r<j$, but then $r>p+q$, and thus $a_{j+k+l,i}=a_{m+q-p,m-r}=0$ and 
$a_{k,i}=a_{m+q,m-r}=0$ by the same reasoning as before. So the first case is 
settled.

\noindent{\bf Second case:}\quad $i=m-r$, $j=m+p$, and $k=m+q$ still with 
$i<j<k$, i.e. $q>p$. These conditions imply $a_{j,k}=0$, and we get
\begin{eqnarray*}
M_{ijk}&=& a_{i,j}a_{i+j+l,k}+a_{k,i}a_{k+i+l,j} \\
&=& a_{m-r,m+p}a_{m-(r-p),m+q}-a_{m-r,m+q}a_{m-(r-q),m+p}.
\end{eqnarray*}
Now we study the relative position of $r$ to $q$: if first $r\geq q$, then
$a_{m-r,m+q}=0$ and $a_{m-r,m+p}=0$. If $r<q$, then $a_{m-(r-q),m+p}=0$ and
following the relative position of $r$ to $p$, either $a_{m-(r-p),m+q}=0$
($r<p$) or $a_{m-r,m+p}=0$ ($r\geq p$). In any case, all terms are zero.\fin

We now come to the last and main point of this section, namely the proof that 
the $m$ family is the only family in weight $l=-m\leq -5$ which satisfies the 
vanishing of all Massey squares which cannot be compensated, i.e. of all 
Massey squares whose vanishing is necessary in order to have an true 
deformation.

Let therefore $\omega$ be a cocycle given by its coefficients $a_{i.j}$. By 
{\bf 2.1}, we have for $i+j\geq m+2$ the usual (or stable) cocycle identity
$a_{i+1,j}+a_{i,j+1}=a_{i,j}$, and for $i+j=m,m+1$ just $a_{i+1,j}+a_{i,j+1}=0$
while there is no equation for lower $i+j$. By {\bf 2.2}, we have $a_{i,j}=0$ 
for $i+j=m+1$, compatible with the foregoing statements. By {\bf 2.0}, the 
coefficients $a_{1,2},\ldots,a_{1,m-1}$ are zero, while by {\bf 2.3}
$a_{1,m+1},a_{1,m+2},\ldots$ may be taken to be zero. Once again, we do not 
impose anything on $a_{1,m}$ in order to use the freedom of choice for a 
coboundary to take the first coefficients (from the diagonal) in the $m+1$st 
column to zero, according to {\bf 2.4}.

The Massey squares which can be 
compensated and thus do not impose conditions on $\omega$ are the $M_{ijk}$
with $i+j=2m+1$, $j+k=2m+1$ or $k+i=2m+1$. 

Let us draw a diagram of the coefficients of $\omega$:   

\vspace{1cm}
\hspace{1cm}
\begin{tabular}{|c||c|c|c|c|c|} \hline
    &   m-2    & m-1  &  m   &  m+1  & m+2      \\  \hline\hline
m-1 & a        &      &      &       &           \\  \hline
m   &  a       &  b   &      &       &            \\  \hline
m+1 & a-b      &  b   &  c   &       &            \\  \hline
m+2 & a-2b     & b-c  &  c   &  0    &             \\  \hline
m+3 & a-3b+c   & b-2c &  c   &  0    &  e           \\  \hline
m+4 & a-4b+3c  & b-3c &  c   &  -e   &  e            \\  \hline
m+5 & a-5b+6c  & b-4c &  c+e & -2e   &               \\  \hline
m+6 &a-6b+10c  &b-5c-e&  c+3e&       &               \\  \hline
m+7 &a-7b+15c+e&b-6c-4e&     &       &               \\  \hline
m+8 &a-8b+21c+5e&      &      &       &               \\  \hline
\end{tabular}
\vspace{1cm}

Let us also expose some Massey squares $M_{ijk}$ (such that no sum of pairs
of indices gives $2m+1$):

\vspace{-1cm}
\hspace{-4cm}
\begin{minipage}[t][5cm][c]{20cm}
\begin{eqnarray*}\scriptstyle
\scriptstyle M_{m-1,m,m+3}&\scriptstyle =&
\scriptstyle a_{m-1,m}a_{m-1,m+3}+a_{m,m+3}a_{m+3,m-1}+
\scriptstyle a_{m+3,m-1}a_{m+2,m}=b(b-2c),    \\ 
\scriptstyle M_{m-1,m+1,m+3}&\scriptstyle=&
\scriptstyle a_{m-1,m+1}a_{m,m+3}+a_{m+1,m+3}a_{m+4,m-1}+
\scriptstyle a_{m+3,m-1}a_{m+2,m+1}=bc,  \\
\scriptstyle M_{m,m+2,m+3}&\scriptstyle=&
\scriptstyle a_{m,m+2}a_{m+2,m+3}+a_{m+2,m+3}
\scriptstyle a_{m+5,m}+a_{m+3,m}a_{m+3,m+1}=e(c-e),    \\ 
\scriptstyle M_{m-2,m,m+2}&\scriptstyle=&
\scriptstyle a_{m-2,m}a_{m-2,m+2}+
\scriptstyle a_{m,m+2}a_{m+2,m-2}=(a-c)(a-2b), \\ 
\scriptstyle M_{m-2,m,m+4}&\scriptstyle=&
\scriptstyle a_{m-2,m}a_{m-2,m+4}+
\scriptstyle a_{m,m+4}a_{m+4,m-2}+a_{m+4,m-2}a_{m+2,m}=a(a-4b+3c),   \\ 
\scriptstyle M_{m-2,m,m+5}&\scriptstyle=&
\scriptstyle a_{m-2,m}a_{m-2,m+5}+a_{m,m+5}a_{m+5,m-2}+
\scriptstyle a_{m+5,m-2}a_{m+3,m}=(a-5b+6c)(a-e).
\end{eqnarray*}
\end{minipage}

\noindent $M_{m-1,m+1,m+3}=bc=0$. Now start a case study:\\

\noindent{\bf First case:} $b=0$, then by $M_{m-2,m,m+2}=0$, either $a=0$ or 
$a=c$. In the first subcase, $M_{m-2,m,m+5}=0$ implies $ce=0$, thus the only
possibly non-zero parameter is $c$ by $M_{m,m+2,m+3}=0$. Note that a non-zero
$c$ corresponds to the $m$ family. In the second subcase, $a=c$ and then 
$M_{m-2,m,m+4}=0$ implies $a=0$. Finally $a=b=c=e=0$ by $M_{m,m+2,m+3}=0$.\\

\noindent{\bf Second case:} $c=0$, then by $M_{m,m+2,m+3}=0$, $e=0$, by
$M_{m-1,m,m+3}=0$, $b=0$, and finally by $M_{m-2,m,m+4}=0$, $a=0$. 

Now it is clear how to perform an induction step showing that the $m$ family
is transmitted to a next two diagonals for example using $M_{m-2,m,m+7}$:
$$M_{m-2,m,m+7}=a_{m-2,m+7}(a_{m-2,m} -a_{m,m+7}+a_{m,m+5}),$$
and by assumption $a_{m-2,m+7}\not= 0$, $a_{m-2,m}=0$, and thus by 
$M_{m-2,m,m+7}=0$, $a_{m,m+7}=a_{m,m+5}=c$ which is the induction step.
This shows that starting from the $(m-2)$nd column, all coeffcients are as in 
the $m$ family. In order to come to lower coefficients, take for example
$$M_{m-3,m+1,m+5}=a_{m-3,m+1}a_{m-2,m+5}+a_{m+1,m+5}a_{m+6,m-3} +
a_{m+5,m-3}a_{m+2,m+1}.$$
Here, $a_{m+1,m+5}=a_{m+2,m+1}=0$ and $a_{m-2,m+5}\not=0$ by assumption, 
therefore $a_{m-3,m+1}=0$ which transmits the $m$-family to the $(m-3)$rd 
column.

Finally, let us argue that the coefficient $a_{1,m}$ cannot be shown to be
zero by the vanishing of Massey squares. Indeed, in order to involve $a_{1,m}$,
the Massey square $M_{ijk}$ (with ordered indices) must have either $i=1$ or
$i+j+l=1$, $j+k+l=1$ or $i+k+l=1$. The first alternative is rather easily seen 
to be zero. Fix $i+j+l=1$ for the second alternative. It describes a 
situation where the coefficient $a_{1,m}$ is multiplied by $a_{i,j}$ with
$i+j=m+1$. This coefficient is zero.  

To summarize, we have the following

\begin{prop} 
The only non-zero cohomology class compatible with the
vanishing of all Massey squares which cannot be compensated, is the
$m$-family
in weight $l=-m\leq -5$, with a possibly non-zero coefficient $a_{1,m}$.
\end{prop} 

\section{Deformations in zero and positive weights} 

\subsection{True deformations in weight $l=0$}

In weight $l\geq 0$, a new phenomenon is happening: we have a relation between 
the Massey squares. Recall that the cocycle coefficients $a_{i,j}$ are 
supposed to be antisymmetric in $i,j$, and that $a_{i,i}$ is set to zero 
for all $i$.

\begin{prop}
Let $i$, $j$, $k$, be three integers, $i,k\geq 2$ and $j\geq 3$. 
We have the relation
$$M_{ijk} + M_{i(j-1)(k+1)} + M_{(i+1)(j-1)k} \,=\, M_{i(j-1)k}.$$
\end{prop}

\pr We have by definition
\begin{eqnarray*}
M_{ijk} + M_{i(j-1)(k+1)} &=& a_{i,j}a_{i+j+l,k} + a_{j,k}a_{j+k+l,i} + \\ 
a_{k,i}a_{i+k+l,j}  + a_{i,j-1}a_{i+j-1+l,k+1} &+& a_{j-1,k+1}a_{j+k+l,i} + 
a_{k+1,i}a_{i+k+1+l,j-1}.
\end{eqnarray*}
 
We transform the terms $a_{j,k}a_{j+k+l,i}+a_{j-1,k+1}a_{j+k+l,i}$, using 
repeatedly the cocycle equation (\ref{***}) to
$$a_{j+k+l,i}(a_{j,k}+a_{j-1,k+1})\,=\,a_{j+k+l,i}a_{j-1,k}\,=\,
a_{j-1,k}(a_{j+k-1+l,i} - a_{j+k-1+l,i+1}).$$ 

We transform the terms $a_{i,j}a_{i+j+l,k} + a_{i,j-1}a_{i+j-1+l,k+1}$, using 
the equations (\ref{***}) to
$$a_{i,j}a_{i+j+l,k} + a_{i,j-1}a_{i+j-1+l,k} - a_{i,j-1}a_{i+j+l,k}.$$

We transform the terms $a_{k,i}a_{i+k+l,j} + a_{k+1,i}a_{i+k+1+l,j-1}$, using 
the equations (\ref{***}) to
$$a_{k,i}a_{i+k+l,j-1} - a_{k,i}a_{i+k+1+l,j-1} + a_{k+1,i}a_{i+k+1+l,j-1}.$$

In these three transformations, the sum of the first term of the first, 
the second term of the second and the first term of the third give together
$$a_{j-1,k}a_{j+k-1+l,i} + a_{i,j-1}a_{i+j-1+l,k} + a_{k,i}a_{i+k+l,j-1}\,=\,
M_{i(j-1)k}.$$

The remaining terms read
$$ - a_{j-1,k}a_{j+k-1+l,i+1} + a_{i,j}a_{i+j+l,k} - a_{i,j-1}a_{i+j+l,k} - 
a_{k,i}a_{i+k+1+l,j-1} + a_{k+1,i}a_{i+k+1+l,j-1}.$$
Here, the second and third term give 
$$a_{i,j}a_{i+j+l,k} - (a_{i,j}a_{i+j+l,k} + a_{i+1,j-1}a_{i+j+l,k}) =
- a_{i+j+l,k}a_{i+1,j-1},$$
while the last two terms give
$$-(a_{k+1,i} + a_{k,i+1})a_{i+k+1+l,j-1} + a_{k+1,i}a_{i+k+1+l,j-1} = 
- a_{i+k+1+l,j-1}a_{k,i+1},$$
still using the equations (\ref{***}).

In summary, the remaining terms give
$$- a_{j-1,k}a_{j+k-1+l,i+1} - a_{i+j+l,k}a_{i+1,j-1} - a_{i+k+1+l,j-1}
a_{k,i+1} = - M_{(i+1)(j-1)k}.$$
This ends the proof of the lemma.\fin

\begin{cor}
$$M_{i(i+1)k} + M_{i(i+2)(k-1)} = M_{i(i+1)(k-1)}.$$
\end{cor}

Observe that also repeated indices may give interesting relations: for example,
for $i=2$, $j=4$ and $k=4$, we get $M_{234}=M_{235}$. It is easily shown by
these relations that the nullity of $M_{23k}$ for all $k\geq 4$ is necessary
for the nullity of all Massey squares, and that the nullity of 
$M_{2rs}$ for all $r,s\geq 3$ is necessary and sufficient for the nullity of
all Massey squares. We believe that the minimal set of Massey squares whose
nullity implies the nullity of all Massey squares is somewhere in between 
these two sets, but we could not get hold on it. 

Now, we will determine all square zero cocycles, i.e. all true 
deformations of $\fm_0$, in weight $l=0$: first of all, the $2$-family 
is such a cocycle. Then, let us suppose that $\omega$ is a non-trivial 
$2$-cocycle which is independent of the $2$-family and has zero Massey squares;
as before, we think of $\omega$ as given by the coefficients $a_{i,j}$, and we 
will distribute letters to its initial terms: $a_{2,3}=a$, $a_{3,4}=b$, and so 
on.

Using equations (\ref{***}), we establish the following diagram which is 
of course valid for all weights $l$; observe that the general expression for
the coefficients in section {\bf 4.5}, proof of proposition $6$, leads for 
general
coefficients $a_{2,3}=:u_2=a, a_{3,4}=:u_3=b, a_{4,5}=:u_4=c$ and so on 
(by linearity) to the formula

\begin{equation}    \label{sum_formula}
a_{i,j}\,=\,\sum_{m=2}^{j-1}(-1)^{m-i}u_m\,
\frac{(j-(m+1))!}{(m-i)!(j-2m+i-1)!},
\end{equation}

which may be used to compute the coefficients in the following diagram 
more easily (than by a recursive formula). 

\vspace{1cm}
\hspace{-1cm}
\begin{tabular}{|c||c|c|c|c|c|c|} \hline
  &  2          &   3     &  4    &  5   & 6 & 7 \\  \hline\hline
3 & a           &         &       &      &   &   \\  \hline
4 &  a          &  b      &       &      &   &   \\  \hline
5 & a-b         &  b      &  c    &      &   &   \\  \hline
6 & a-2b        & b-c     &  c    &  d   &   &   \\  \hline
7 & a-3b+c      & b-2c    &  c-d  &  d   & e &   \\  \hline
8 & a-4b+3c     & b-3c+d  &  c-2d &  d-e & e & f \\  \hline
9 & a-5b+6c-d   & b-4c+3d &c-3d+e &  d-2e&e-f& f \\  \hline
10 &a-6b+10c-4d &b-5c+6d-e&c-4d+3e&d-3e+f&e-2f&  \\  \hline
11&a-7b+15c-10d+e&b-6c+10d-4e&c-5d+6e-f&d-4e+3f&&\\  \hline
12&a-8b+21c-20d+5e&b-7c+15d-10e+f&c-6d+10e-4f&&& \\  \hline
13&a-9b+28c-35d+15e-f&b-8c+21d-20e+5f&&&&        \\  \hline   
14&a-10b+36c-56d+35e-6f& &&&&                    \\  \hline  
\end{tabular}
\vspace{1cm}

>From now on, we consider weight $l=0$.

Order the Massey squares by their {\it level}, i.e. we say that 
$M_{ijk}$ has level $i+j+k$. Computing Massey squares and setting them 
equal to zero gives an infinite family of homogeneous quadratic equations 
for the infinite family of variables $a,b,c,d,\ldots$.

In level $9$, the only Massey square is $M_{234}$, and its nullity gives
$$3b^2-bc - 2ac = 0.$$
In level $10$, the only Massey square is $M_{235}$, and its nullity gives
the same equation.
In level $11$, there are Massey squares $M_{236}$ and $M_{245}$, and their
nullity gives (possibly by subtracting the previous equation) in both cases
$$2ad - 4bc -bd + 6c^2 -cd = 0.$$ 
In level $12$, there are Massey squares $M_{237}$, $M_{246}$ and $M_{345}$,
and their nullity gives (possibly by subtracting the previous equations) 
in all cases
$$-3bd + 4c^2 - 3cd = 0.$$

Going higher in this hierarchy of equations and variables, there are at each
new level some (possibly) linear independent equation. Proposition $9$ only
tells us that the nullity of $M_{2rs}$ with $2<r<s$ is enough in order to have
all Massey squares zero. We don't know which of these equations are in fact
the independent one's. 

In Massey square level $14$, we arrive at $5$ 
equations for the $5$ variables $a,b,c,d,e$, which read (after subtracting
at each step multiples of the previous equations):

\begin{eqnarray*}
3b^2-bc - 2ac &=& 0  \\
2ad - 4bc -bd + 6c^2 -cd &=& 0\\
-3bd + 4c^2 - 3cd &=& 0 \\
e(-2a+3b-d)+5bd-15cd+10d^2 &=& 0 \\
e(-6a+15b-4c-11d)-55cd+50d^2+15bd &=& 0
\end{eqnarray*}

The discussion of these equations (either by hand or by a system computing a 
Gr\"obner basis for the homogeneous polynomials) gives as non-zero solutions 
the $2$-family and one other family with coefficients $a=\frac{1}{6}$,
$b=\frac{1}{60}$, $c=\frac{1}{420}$, etc. We describe this family from 
another point of view in subsection {\bf 5.3}, which will show that this family
must verify all equations and not only the five equations we wrote down. 
These are the only square zero solutions in weight $0$, and we have determined 
all true deformations in this case.

\subsection{True deformations in weight $l>0$}

In the weight $l=1$ case, we get from the same diagram as in weight $l=0$
up to Massey level $15$ (where we took only the equations of type $M_{23k}$ in 
order to simplify) six homogeneous quadratic equations in six variables
which read:

\begin{eqnarray*}
\scriptstyle -3ac+4b^2-3bc &\scriptstyle=&\scriptstyle 0 \\
\scriptstyle -5bc-2bd+10c^2-4cd+3ad &\scriptstyle=&\scriptstyle 0 \\
\scriptstyle 5c^2-4bd+2be+ec-6cd &\scriptstyle=&\scriptstyle 0 \\
\scriptstyle e(-3a+11b-5d+3c)-6bd+15c^2-39cd+20d^2 &
\scriptstyle=&\scriptstyle 0 \\
\scriptstyle e(-9a+35b-35d)+f(-4b+2c+2d)-6bd+30c^2-111cd+90d^2 &
\scriptstyle=&\scriptstyle 0 \\
\scriptstyle e(-18a+75b+11c-186d+35e-6f)+f(3a+20c-24b+16d)-4bd+
50c^2-234cd+252d^2 &\scriptstyle=&\scriptstyle 0
\end{eqnarray*}

By a computation with MUPaD which determines a Gr\"obner basis for the 
homogeneous polynomials corresponding to these equations (actually we took 
here all equations of type $M_{2rs}$), one obtains
as (non-zero) solutions the $2$-family, a solution $a=0, b=0, c=0, d=0$ and
$e=1$, $f=\frac{35}{6}$, and a further solution $a=1$, $b=\frac{1}{7}$,
$c=\frac{1}{42}$, $d=\frac{1}{231}$, $e=\frac{5}{21\cdot 286}$, 
$f=\frac{1}{21\cdot 286}$. The solution with $e=1$ and $f=\frac{35}{6}$
does not survive the next level of Massey squares. 

But the solution starting with $a=1$, $b=\frac{1}{7}$ continues with
$g=\frac{1}{29172}$, $h=\frac{1}{138567}$, $i=\frac{1}{646646}$, 
$j=\frac{5}{14872858}$, $k=\frac{1}{13520780}$. We will describe this family
from a different point of view in section {\bf 5.3}, and we will show there
(implicitly) that this solution survives to infinity.

\begin{prop}
In weight $l=1$, there are exactly two non-equivalent true deformations.
\end{prop}

The problem of determining the explicit square zero cocycles in each weight 
$l$ case seems to be a rich problem. We tried to say something about the rank
of the finite Jacobimatrix either associated to the set of equations of type 
$M_{2rs}=0$ for all $r,s$, or to those of type $M_{23r}$ for all $r$, when
we truncate the number of variables and consider only those equations 
involving these variables. With this matrix, it is obvious that the set of 
solutions (of the truncated problem) is an algebraic variety of dimension
greater or equal to $1$ (because the equations are homogeneous), but we
couldn't decide whether the dimension drops down to $1$ in each weight. 
In fact, within the possibilities of our computer, we computed (using all 
equations of type $M_{2rs}=0$) the dimension of this
variety as far as possible for $l=2$ and it remained $2$. Is the set of 
solutions always a variety ? Is it always of finite dimension ? Can we give
asymptotics or bounds or a formula for the dimension ?

\subsection{Identifying the cocycles and the deformed algebras in weight zero}

We now construct deformations from the previously determined weight $l$
$2$-cocycles given by their coefficients $a_{i,j}$ in the following way: 
using still the $e_i$ for $i\geq 1$ as a basis, the deformed bracket reads
$$[e_i,e_j]_t\,=\,[e_i,e_j]\,+\,ta_{i,j}e_{i+j+l}.$$
It is clear that all square zero $2$-cocycles give in this way true 
deformations of $\fm_0$ for which only the linear term is (possibly) non-zero; 
this means in particular that the bracket $[-,-]_t$ satisfies the Jacobi 
identity without adding terms containing higher powers in $t$.

The weight $l=0$ case is the most interesting, because here deformations give
automatically rise to $\N$-graded Lie algebras which must fit in the 
classification \cite{Fia1}. In this classification, the three $\N$-graded
Lie algebras where $e_1$ has non-zero brackets with all other basis 
elements are
\begin{itemize}
\item[(1)] $\fm_0$; brackets: $[e_1,e_i]=e_{i+1}$ for all $i\geq 2$,
\item[(2)] $\fm_2$; brackets: $[e_1,e_i]=e_{i+1}$ for all $i\geq 2$, 
$[e_2,e_j]=e_{j+1}$ for all $j\geq 3$,
\item[(3)] $L_1$; brackets: $[e_i,e_j]=(j-i)e_{i+j}$ for all $i,j\geq 1$.
\end{itemize}
The complete set of infinitesimal deformations of $L_1$ and the complete set
of formal deformations of $L_1$ is given in \cite{Fia2}. Let us consider in
this section the same problem for $\fm_0$ in weight $l=0$.

Taking as $2$-cocycle the $2$-family, we get in weight $l=0$ a Lie algebra
$\fm_0^1(t)$ which must be $\N$-graded and which is easily seen to be generated
by $e_1$ and $e_2$: indeed, $[e_1,e_i]_t=[e_1,e_i]$ for all $i\geq 2$. The 
complete relations for $\fm_0^1(t)$ are
$$\left\{\begin{array}{ccccc} [e_1,e_i]_t&=&e_{i+1} &\forall& i\geq 2 
\\ \,[e_2,e_j]_t &=& te_{j+2} &\forall& j\geq 3 
\end{array}\right.$$
Thus this family describes the deformation of $\fm_0$ to $\fm_2$.

Now there is also a cocycle describing the deformation of $\fm_0$ to $L_1$:
observe that the generators $e_i$ for $i\geq 1$ of $\fm_0$ do not
satisfy the right relations, seen as elements of $L_1$. Therefore, one must
first perform a change of base: let us define $\tilde{e}_1=e_1$, 
$\tilde{e}_2=e_2$, $\tilde{e}_3=e_3$, $\tilde{e}_4=\frac{1}{2}e_4$ and in 
general
$$\tilde{e}_i\,=\,\frac{1}{(i-2)!}e_i$$
for all $i\geq 5$.
Then the relations are easily computed to be
$$[\tilde{e}_1,\tilde{e}_i]=(i-1)\tilde{e}_{i+1}.$$
The cocycle relation (\ref{***}) for a cocycle given by coefficients
$b_{i,j}$ transforms in the new basis to
$$(j-1)b_{j+1,k}\,+\,(k-1)b_{j,k+1}\,=\,(j+k+l-1)b_{j,k}.$$
It is easy to check that the $2$-cochain given by the coefficients
$b_{i,j}=(j-i)$ is indeed a $2$-cocycle for $l=0$. One also easily checks that
the $2$-cocycle $b_{i,j}=(j-i)$ is of Massey square zero, i.e.
$$b_{i,j}b_{i+j,k} + b_{j,k}b_{j+k,i} + b_{k,i}b_{k+i,j}=0.$$
It therefore determines a deformation $\fm_0^2(t)$ of $\fm_0$ in weight $l=0$ 
to $L_1$, and we showed in the previous section that these are all possible 
deformations of $\fm_0$ in weight $0$. 

Let us finish with the identification of the deformations in weight $l=1$. 
It is clear that the $2$-family leads to a non-trivial true deformation. This
is then a weight $1$ variant of the Lie algebra $\fm_2$. The other cocycle, 
determined using MUPaD, is more interesting. 
Indeed, there is a general procedure 
of constructing positive weight, true deformations for $\fm_0$: consider
the Lie algebra $L_1$, with its generators $e_1$, $e_2$, $e_3$, etc and its 
relations $[e_i,e_j]=(j-i)e_{i+j}$ for all $i,j\geq 1$. Define a Lie algebra 
$L_1\{2\}$ by generators $e_1$, $e_3$, $e_4$ etc (the suppression of $e_2$
is indicated by $\{2\}$ in the notation !) and the relations of $L_1$ for the
remaining generators. Introduce a new basis $f_1:=e_1$, $f_2:=e_3$, 
$f_3:=e_4$ etc, and another new basis $g_1:=e_1$, $g_k:=(k-1)!\,\,f_k$ for all
$k>1$. We compute the relations to
$$[g_1,g_k]\,=\,g_{k+1},\,\,\,\,\,\,\,\,[g_k,g_{k+1}]\,=\,
\frac{k!(k-1)!}{(2k+1)!}g_{2k+2}.$$  
When interpreted as a deformation of $\fm_0$, one can then compute the 
coefficients $a_{i,j}$ of the corresponding $2$-cocycle. One obtains
$a_{2,3}=\frac{2!1!}{5!}=\frac{1}{60}\cdot 1$, $a_{3,4}=\frac{3!2!}{7!}=
\frac{1}{60}\cdot \frac{1}{7}$, $a_{4,5}=\frac{4!3!}{9!}=\frac{1}{60}\cdot 
\frac{1}{42}$, $a_{5,6}=\frac{5!4!}{11!}=\frac{1}{60}\cdot\frac{1}{231}$,
$a_{6,7}=\frac{6!5!}{13!}=\frac{1}{60}\cdot\frac{5}{21\cdot 286}$,
$a_{7,8}=\frac{7!6!}{15!}=\frac{1}{60}\cdot\frac{1}{21\cdot 286}$, and so on.
Thus, this deformation has as its infinitesimal cocycle the cocycle we 
determined using MUPaD before, up to a factor $\frac{1}{60}$.
By construction, it is clear that it defines
a cocycle and a true deformation, because $L_1\{2\}$ is a Lie algebra.

In the same way, one can define $L_1\{m\}$ by the span of the vectors
$e_1$, $e_{m+1}$, $e_{m+2}$, etc for any $m>2$, in other words, by the 
suppression of all basis vectors from $e_2$ up to and including $e_{m+1}$.
This gives a non-trivial 
true deformation of weight $m-1$. Together with the deformation given by the 
$2$-family, these two constitute two independent true deformations in any 
positive weight.

\end{document}